\documentclass[bsl]{asl}

\title[Relations between cardinals in the absence of AC]
{\large Relations between some cardinals\\[0.7ex] in the absence of
the Axiom of Choice\\[4ex] \scriptsize\it{Dedicated to the memory of
Prof.\;Hans L\"auchli}}


\author[Lorenz Halbeisen]{Lorenz Halbeisen${}^1$}
\revauthor{Halbeisen, Lorenz}

\address{Department of Mathematics\\
Universtity of California at Berkeley\\ Berkeley, CA 94720, USA}

\email{halbeis@math.berkeley.edu}

\thanks{${}^1$I would like to thank the {\it Centre de
Recerca Matem\`atica (Barcelona)\/} and the {\it Swiss National
Science Foundation\/} for their support.}

\author[Saharon Shelah]{Saharon Shelah${}^2$}
\revauthor{Shelah, Saharon}

\address{Institute of Mathematics\\
The Hebrew University Jerusalem\\ Jerusalem 91904, Israel}

\email{shelah@math.huji.ac.il}

\thanks{${}^2$Partially supported by the {\it Israel Basic
Research Fund}, founded by the Israel Academy; Publ.\,No.\,699}

\newcommand \mdom {{\mbox{dom}}}
\newcommand \la {{\langle}}
\newcommand \ra {{\rangle}}
\newcommand \subs {{\subseteq}}
\newcommand \supp {\operatorname{supp}}
\newcommand \lesseq {\preccurlyeq}

\newcommand \cP {{\mathcal{P}}}
\newcommand \cF {{\mathcal{F}}}
\newcommand \cC {{\mathcal{C}}}
\newcommand \cG {{\mathcal{G}}}
\newcommand \cV {{\mathcal{V}}}
\newcommand \fm {{\mathfrak{m}}}
\newcommand \fn {{\mathfrak{n}}}

\newcommand \fp {{\mathfrak{p}}}
\newcommand \fq {{\mathfrak{q}}}
\newcommand \fzw {{\mathfrak{2}}}
\newcommand \seq {\operatorname{seq}}
\newcommand \Part {\operatorname{part}}
\newcommand \Sym {\operatorname{sym}}
\newcommand \Fin {\operatorname{fin}}
\newcommand \Fix {\operatorname{fix}}
\newcommand \On {\operatorname{On}}

\newcommand \oneone {{\hspace{-0.1ex}}{\rule[0.3ex]{0.8mm}{0.1mm}}
                                       {\hspace{-0.08ex}}}
\newcommand \Seq {{\operatorname{seq}}^{\scriptscriptstyle{1}\oneone
                   \scriptscriptstyle{1}}}

\newcommand \incomp {{\parallel}}
\newcommand \N {{\mathbb{N}}}

\theoremstyle{definition}
\newtheorem {nummer}{ }[section]

\newtheorem {prop}[nummer]{\sc{Proposition}}

\newtheorem {fct}[nummer]{\sc{Fact}}

\newtheorem {subnummer}{ }[subsection]

\newtheorem {sprop}[subnummer]{\sc{Proposition}}
\newtheorem {slm}[subnummer]{\sc{Lemma}}
\newtheorem {sfct}[subnummer]{\sc{Fact}}

\newcommand \rmk {{\sc{Remark:}}\hspace*{3mm}}
\newcommand \defi {{\sc{Definition:}}\hspace*{3mm}}
\newcommand \notation {{\sc{Notation:}}\hspace*{3mm}}

\begin{document}

\begin{abstract}
If we assume the axiom of choice, then every two cardinal numbers are
comparable. In the absence of the axiom of choice, this is no longer
so. For a few cardinalities related to an arbitrary infinite set, we
will give all the possible relationships between them, where possible
means that the relationship is consistent with the axioms of set
theory. Further we investigate the relationships between some other
cardinal numbers in specific permutation models and give some results
provable without using the axiom of choice.
\end{abstract}

\maketitle

\section{Introduction}

Using the axiom of choice, Felix Hausdorff proved in 1914 that there
exists a partition of the sphere into four parts, $S=A\,\dot\cup\,
B\,\dot\cup\,C\,\dot\cup\,E$, such that $E$ has Lebesgue measure $0$,
the sets $A,\,B,\,C$ are pairwise congruent and $A$ is congruent to
$B\,\dot\cup\,C$ ({\sl cf.}\;\cite{Hausdorff1} or \cite{Hausdorff2}).
This theorem later became known as Hausdorff's paradox. If we want to
avoid this paradox, we only have to reject the axiom of choice. But
if we do so, we will run into other paradoxical situations. For
example, without the aid of any form of infinite choice we cannot
prove that a partition of a given set $m$ has at most as many parts
as $m$ has elements. Moreover, it is consistent with set theory that
the real line can be partitioned into a family of cardinality
strictly bigger than the cardinality of the real numbers (see
Fact\;\ref{fct:partition}).

Set theory without the axiom of choice has a long tradition and a lot
of work was done by the Warsaw School between 1918 and 1940.
Although, in 1938, Kurt G\"odel proved in \cite{Godel} the
consistency of the axiom of choice with the other axioms of set
theory, it is still interesting to investigate which results remain
unprovable without using the axiom of choice ({\sl
cf.}\;\cite{Lauchli2}).

In 1963, Paul Cohen proved with his famous and sophisticated forcing
technique, that it is also consistent with the other axioms of set
theory that the axiom of choice fails ({\sl cf.}\;\cite{Cohen}). Also
with a forcing construction, Thomas~Jech and Anton{\'\i}n~Sochor
could show in \cite{JechSochor} that one can embed the permutation
models (these are models of set theory with atoms) into well-founded
models of set theory. So, to prove consistency results in set theory,
it is enough to build a suitable permutation model.

We will investigate the relationships between some infinite cardinal
numbers. For four cardinal numbers---which are related to an
arbitrary given one---we will give all the possible relationships
between two of them; where possible means that there exists a model
of set theory in which the relationship holds.  For example it is
possible that there exists an infinite set $m$ such that the
cardinality of the set of all finite sequences of $m$ is strictly
smaller than the cardinality of the set of all finite subsets of $m$.
On the other hand, it is also possible that there exists an infinite
set $m'$ such that the cardinality of the set of all finite sequences
of $m'$ is strictly bigger than the cardinality of the power-set of
$m'$. In a few specific permutation models, like the basic Fraenkel
model and the ordered Mostowski model, we will investigate also the
relationships between some other cardinal numbers. Further we give
some results provable without using the axiom of choice and show that
some relations imply the axiom of choice.

\section{Definitions, notations and basic facts}

First we want to define the notion of a cardinal number and for this
we have to give first the definition of ordinal numbers.

\defi A set $\alpha$ is an {\bf ordinal} if and only if every
element of $\alpha$ is a subset of $\alpha$ and $\alpha$ is
well-ordered by $\in$.

Now let $V$ be a model for ZF (this is Zermelo-Fraenkel's set theory
without the axiom of choice) and let $\On:=\{\alpha\in
V:\alpha\;\text{is an ordinal}\}$; then $\On$ is a proper class in
$V$. It is easy to see that if $\alpha\in\On$, then also $\alpha
+1:=\alpha\cup\{\alpha\}\in\On$. An ordinal $\alpha$ is called a {\bf
successor ordinal} if there exists an ordinal $\beta$ such that
$\alpha =\beta +1$ and it is called a {\bf limit ordinal} if it is
neither a successor ordinal nor the empty-set.

By transfinite recursion on $\alpha\in\On$ we can define $V_\alpha$
as follows: $V_\emptyset :=\emptyset$, $V_{\alpha +1}=\cP(V_\alpha)$
and $V_\alpha:=\bigcup_{\beta\in\alpha}V_\beta$ when $\alpha$ is a
limit ordinal. Note that by the axiom of power-set and the axiom of
replacement, for each $\alpha\in\On$, $V_\alpha$ is a set in $V$. By
the axiom of foundation we further get $V:=\bigcup_{\alpha\in
\On}V_\alpha$ ({\sl cf.}\;\cite[Theorem\;4.1]{Kunen}).

Let $m$ be a set in $V$, where $V$ is a model of ZF, and let
$\mathfrak{C}(m)$ denote the {\bf class of all sets $x$}, such that
there exists a {\bf one-to-one mapping from $x$ onto} $m$. We define
the cardinality of $m$ as follows.

\defi For a set $m$, let $\fm:=\mathfrak{C}(m)\cap V_\alpha$, where
$\alpha$ is the smallest ordinal such that
$V_\alpha\cap\mathfrak{C}(m)\neq\emptyset$. The set $\fm$ is called
the {\bf cardinality} of $m$ and a set $\fn$ is called a {\bf
cardinal number} (or simply a {\bf cardinal}) if it is the
cardinality of some set.

Note that a cardinal number is defined as a set.

A cardinal number $\fm$ is an {\bf aleph} if it contains a
well-ordered set. So, the cardinality of each ordinal is an aleph.
Remember that the axiom of choice is equivalent to the statement that
each set can be well-ordered. Hence, in ZFC (this is
Zermelo-Fraenkel's set theory with the axiom of choice), every
cardinal is an aleph; and vice versa, if every cardinal is an aleph,
then the axiom of choice holds.

If we have a model $V$ of ZF in which the axiom of choice fails, then
we have more cardinals in $V$ than in a model $M$ of ZFC. This is
because all the ordinals are in $V$ and, hence, the alephs as well.

\notation We will use fraktur-letters to denote cardinals and
$\aleph$'s to denote the alephs. For finite sets $m$, we also use
$|m|$ to denote the cardinality of $m$. Let $\N
:=\{0,1,2,\ldots\}$ be the set of all natural numbers and let
$\aleph_0$ denote its cardinality. We can consider $\N$ also as the
set of finite ordinal numbers, where $n=\{0,1,\ldots,n-1\}$ and
$0=\emptyset$. For a natural number $n\in\N$, we will not distinguish
between $n$ as an ordinal number and the cardinality of $n$. Further,
the ordinal number $\omega$ denotes the order-type (with respect to
$<$) of the set $\N$.

Now we define the order-relation between cardinals.

\defi We say that the cardinal number $\fp$
is {\bf less than or equal} to the cardinal number $\fq$ if and only
if for any $x \in\fp$ and $y \in\fq$ there is a {\bf one-to-one
mapping from $x$ into $y$}.

\notation If $\fp$ is less than or equal to the cardinal number
$\fq$, we write $\fp\leq\fq$. We write $\fp<\fq$ for $\fp\leq\fq$ and
$\fp\neq\fq$. If neither $\fp\leq\fq$ nor $\fq\leq\fp$ holds, then we
say that $\fp$ and $\fq$ are {\bf incomparable} and write
$\fp\incomp\fq$. For $x \in\fp$ and $y \in\fq$ we write: $x\lesseq y$
if $\fp\le\fq$ and $x\not\lesseq y$ if $\fp\nleq\fq$ ({\sl cf.}\;also
\cite[p.\,27]{Kunen}). Notice that $x\lesseq y$ iff there exists a
one-to-one function from $x$ into $y$.

Another order-relation which we will use at a few places and which
was first introduced by Alfred~Tarski ({\sl
cf.}\;\cite{LindenbaumTarski}) is the following.

\defi For two cardinal numbers $\fp$ and $\fq$ we write $\fp\leq^*\fq$
if there are non-empty sets $x \in\fp$ and $y \in\fq$ and a function
from $y$ {\bf onto} $x$.

Notice, that for infinite cardinals $\fp$ and $\fq$, we must use the
axiom of choice to prove that $\fp\leq^*\fq$ implies $\fp\leq\fq$
(see {\sl e.g.}\;\cite{LoriShelah}). In general, if we work in ZF,
there are many relations between cardinals which do not exist if we
assume the axiom of choice ({\sl cf.}\;\cite{LoriShelah}); and
non-trivial relations between cardinals become trivial with the axiom
of choice (see also \cite{Lauchli1} or \cite{Sierpinski}).

The main tool in ZF to show that two cardinals are equal is the

{\sc{Cantor-Bernstein Theorem:}} {If $\fp$ and $\fq$ are cardinals
with $\fp\leq\fq$ and $\fq\leq\fp$, then $\fp=\fq$.}

(For a proof see\;\cite{Jechbook} or \cite{Bachmann}.)

Notice that for $x\in\fp$ and $y\in\fq$ we have $x\lesseq
y\not\lesseq x$ is equivalent to $\fp <\fq$, and if $x\lesseq
y\lesseq x$, then there exists a one-to-one mapping from $x$ onto
$y$.

A result which gives the connection between the cardinal numbers and
the $\aleph$'s is

{\sc{Hartogs Theorem:}} {For every cardinal number $\fm$, there
exists a least aleph, denoted by $\aleph(\fm)$, such that
$\aleph(\fm) \nleq \fm$.}

(This was proved by Friedrich~Hartogs in \cite{Hartogs}, but a proof
can also be found in~\cite{Jechbook} or in~\cite{Bachmann}.)

Now we will define ``infinity''.

\defi A cardinal number is called {\bf finite} if it is the cardinality
of a natural number, and it is called {\bf infinite} if it is not
finite.

There are some other degrees of infinity ({\sl cf.}\;{\sl
e.g.}\;\cite{Goldstern} or \cite{SpisiakVojtas}), but we will use
only ``infinite'' for ``not finite'' and as we will see, most of the
infinite sets we will consider in the sequel will be Dedekind finite,
where a cardinal number $\fm$ is called {\bf Dedekind finite} if
$\aleph_0\nleq\fm$.

There are also many weaker forms of the axiom of choice (we refer the
reader to \cite{HowardRubin}). Concerning the notion of Dedekind
finite we wish to mention five related statements.
$$\begin{array}{rl} \text{AC:} & \text{``The Axiom of Choice'';}\\
\fzw\fm=\fm: & \text{``For every infinite cardinal $\fm$ we have
$\fzw\fm=\fm$'';}\\ C(\aleph_0,\infty): & \text{``Every countable
family of non-empty sets has a}\\ & \text{\ \;choice function'';}\\
C(\aleph_0,<\aleph_0): & \text{``Every countable family of non-empty
finite sets}\\ & \text{\ \;has a choice function'';}\\ W_{\aleph_0}:
& \text{``Every Dedekind finite set is finite''.}\end{array}$$

We have the following relations (for the references
see~\cite{HowardRubin}):
$$\text{AC}\,\Rightarrow\,\fzw\fm=\fm\,\Rightarrow\,W_{\aleph_0}
\,\Rightarrow\,C(\aleph_0,<\aleph_0)\ \;\text{and}\ \;
\text{AC}\,\Rightarrow\,C(\aleph_0,\infty)
\,\Rightarrow\,W_{\aleph_0},$$ but on the other hand have
$$\text{AC}\,\nLeftarrow\,\fzw\fm=\fm\,\nLeftarrow\,W_{\aleph_0}
\,\nLeftarrow\,C(\aleph_0,<\aleph_0)\ \;\text{and}\ \;
\text{AC}\,\nLeftarrow\,C(\aleph_0,\infty)
\,\nLeftarrow\,W_{\aleph_0},$$ and further $\fzw\fm=\fm\;
\nRightarrow\;C(\aleph_0,\infty)\;\nRightarrow\;\fzw\fm=\fm$.

\section{Cardinals related to a given one}\label{sec:cardinals}

Let $m$ be an arbitrary set and let $\fm$ denote the cardinality of
$m$. In the following we will define some cardinalities which are
related to the cardinal number $\fm$.

Let $[m]^2$ be the set of all $2$-element subsets of $m$ and let
$[\fm]^2$ denote the cardinality of the set $[m]^2$.

Let $\Fin (m)$ denote the set of all finite subsets of $m$ and
let $\Fin (\fm)$ denote the cardinality of the set $\Fin (m)$.

For a natural number $n$, $\Fin(m)^n$ denotes the set $\{\la
e_0,\ldots,e_{n-1}\ra : \forall i<n (e_i\in \Fin(m))\}$ and $\Fin
(\fm)^n$ denotes its cardinality.

For a natural number $n$, $\Fin^{n+1} (m)$ denotes the set
$\Fin (\Fin^n (m))$, where $\Fin^0 (m):=m$,
and $\Fin^{n+1} (\fm)$ denotes its cardinality.

Let $m^2:=m\times m=\{\la x_1,x_2\ra: \forall i<2(x_i\in m)\}$ and
let $\fm^2=\fm\cdot\fm$ denote the cardinality of the set $m^2$.

Let $\Seq (m)$ denote the set of all finite one-to-one sequences of
$m$, which is the set of all finite sequences of elements of $m$ in
which every element appears at most once, and let $\Seq (\fm)$ denote
the cardinality of the set $\Seq (m)$.

Let $\seq (m)$ denote the set of all finite sequences of $m$
and let $\seq (\fm)$ denote the cardinality
of the set $\seq (m)$.

Finally, let $\cP (m)$ denote the power-set of $m$, which is the set
of all subsets of $m$, and let $\fzw^\fm$
denote the cardinality of $\cP (m)$.

In the sequel, we will investigate the relationships between these
cardinal numbers.

\section{Cardinal relations which imply the axiom of choice}
\label{sec:choice}

First we give some cardinal relations which are well-known to be
equivalent to the axiom of choice. Then we show that also a weakening
of one of these relations implies the axiom of choice.

The following equivalences are proved by Tarski in 1924. For the
historical background we refer the reader to~\cite[4.3]{Moore}.

\begin{prop}\label{prop:tarski}
The following conditions are equivalent to the axiom of
choice:\hfill\\ \hspace*{3ex}(1) $\fm\cdot\fn=\fm+\fn$ for every
infinite cardinal $\fm$ and $\fn$\hfill\\ \hspace*{3ex}(2)
$\fm=\fm^2$ for every infinite cardinal $\fm$\hfill\\
\hspace*{3ex}(3) If $\fm^2=\fn^2$, then $\fm=\fn$\hfill\\
\hspace*{3ex}(4) If $\fm<\fn$ and $\fp<\fq$, then
$\fm+\fp<\fn+\fq$\hfill\\ \hspace*{3ex}(5) If $\fm<\fn$ and
$\fp<\fq$, then $\fm\cdot\fp<\fn\cdot\fq$\hfill\\ \hspace*{3ex}(6) If
$\fm+\fp<\fn+\fp$, then $\fm<\fn$\hfill\\ \hspace*{3ex}(7) If
$\fm\cdot\fp<\fn\cdot\fp$, then $\fm<\fn$
\end{prop}

(The proofs can be found in \cite{Tarski}, in \cite{Bachmann}
or in \cite{Sierpinski}.)

As a matter of fact we wish to mention that Tarski observed that the
statement $$\text{If\ }2\fm<\fm+\fn,\text{\ then\ }\fm<\fn$$ is
equivalent to the axiom of choice, while the proposition: $$\text{If\
}2\fm>\fm+\fn,\text{\ then\ }\fm>\fn$$ can be proved without the aid
of the axiom of choice ({\sl cf.}\;\cite[p.\,421]{Sierpinski}).

To these cardinal equivalences mentioned above, we will now add two more:

\begin{prop}\label{prop:lori}
The following conditions are equivalent to the axiom of
choice:\hfill\\ \hspace*{3ex}(1) For every infinite cardinal $\fm$ we
have $[\fm]^2=\fm$\hfill\\ \hspace*{3ex}(2) For every infinite
cardinal $\fm$ we have $[\fm]^2=\fm$ or $\fm^2=\fm$
\end{prop}

\begin{proof} The proof is essentially the same as Tarski's proof that
the axiom of choice follows if $\fm^2=\fm$ for all infinite cardinals
$\fm$ ({\sl cf.}\;\cite{Tarski}).

Tarski proved in~\cite{Tarski} ({\sl cf.}\;also \cite{Sierpinski})
the following relation for infinite cardinals $\fm$: $$\fm+\aleph
(\fm)=\fm\cdot\aleph (\fm) \text{\ {\it implies\/}\ }\fm <\aleph
(\fm).$$ Notice that $\fm
<\aleph (\fm)$ implies that every set $m\in\fm$ can be well-ordered.
Therefore it is sufficient to show that (2), which is weaker than
(1), implies that for every infinite cardinal number $\fm$ we have
$\fm
<\aleph (\fm)$.

First we show that for two infinite cardinal numbers $\fm$ and $\fn$
we have $\fm+\fn\le\fm\cdot\fn$. For this, let $\fm_1$ and $\fn_1$ be
such that $\fm=\fm_1 +1$ and $\fn=\fn_1 +1$. Now we get
$$\fm\cdot\fn=(\fm_1 +1)\cdot(\fn_1 +1)=\fm_1\cdot\fn_1 +\fm_1 +\fn_1
+1 \ge 1+\fm_1 +\fn_1 +1=\fm +\fn.$$ It is easy to compute, that
$$[\fm +\aleph (\fm)]^2 = [\fm]^2 +\fm\aleph (\fm)+[\aleph
(\fm)]^2,$$ and $$(\fm +\aleph (\fm))^2 = \fm^2 +2\fm\aleph
(\fm)+\aleph (\fm)^2.$$ Now we apply the assumption (2) to the
cardinal $\fm+\aleph (\fm)$. If $[\fm +\aleph (\fm)]^2 =\fm +\aleph
(\fm)$, we get $\fm\aleph (\fm)\le \fm +\aleph (\fm)$ which implies
(by the above, according to the Cantor-Bernstein Theorem) $\fm\aleph
(\fm)= \fm +\aleph (\fm)$. By the result of Tarski mentioned above we
get $\fm < \aleph (\fm)$. The case when $(\fm +\aleph (\fm))^2 = \fm
+\aleph (\fm)$ is similar. So, if the assumption (2) holds, then we
get $\fm <\aleph (\fm)$ for every cardinal number $\fm$ and
therefore, each set $m$ can be well-ordered, which is equivalent to
the axiom of choice. \end{proof}

\section{A few relations provable in ZF}\label{sec:provable}

In this section we give some relationships between the cardinal
numbers defined in section\;\ref{sec:cardinals} which are provable
without using the axiom of choice.

The most famous one is the

{\sc{Cantor Theorem:}} {For any cardinal number $\fm$ we have $\fm
<\fzw^\fm$.}

(This is proved by Georg~Cantor in \cite{Cantor}, but a proof can
also be found in~\cite{Jechbook} or \cite{Bachmann}.)

Concerning the relationship between ``$\le^*$'' and ``$\le$'', it is
obvious that $\fp\le\fq$ implies $\fp\le^*\fq$. The following fact
gives a slightly more interesting relationship.

\begin{fct}\label{fct:*}
For two arbitrary cardinals $\fn$ and $\fm$ we have
$\fn\le^*\fm\rightarrow \fzw^\fn\le\fzw^\fm$.
\end{fct}

(For a proof see {\sl e.g.}\;\cite{Sierpinski} or \cite{Bachmann}.)

The following two facts give a list of a few obvious relationships.

\begin{fct}For every cardinal $\fm$ we have:\hfill\\
\hspace*{3ex}(1) $\fm^2\le\Fin^2 (\fm)$\hfill\\ \hspace*{3ex}(2)
$\Seq (\fm) \le\Fin^2 (\fm)$\hfill\\ \hspace*{3ex}(3) $\Seq
(\fm)\le\seq (\fm)$\hfill\\ \hspace*{3ex}(4) If $\fm$ is infinite,
then $\fzw^{\aleph_0}\leq\fzw^{\Fin (\fm)}$
\end{fct}

\begin{proof} First take an arbitrary set $m\in\fm$. For (1) note that a set
$\la x_1,x_2\ra\in m^2$ corresponds to the set
$\{\{x_1\},\{x_1,x_2\}\}\in\Fin^2 (m)$. For (2) note that a finite
one-to-one sequence $\la a_0,a_1,\ldots,a_n\ra$ of $m$ can always be
written as
$\{\{a_0\},\{a_0,a_1\},\ldots,\{a_0,\ldots,\linebreak[2]a_n\}\}$,
which is an element of $\Fin^2 (m)$. The relation (3) is trivial. For
(4) let $E_n:=\{e\subs m: |e|=n\}$, where $n\in\N$. Because $m$ is
assumed to be infinite, every $x\subs\N$ corresponds to a set
$F_x\in\cP(\Fin (m))$ defined by $F_x:=\bigcup\{E_n: n\in x\}$.
\end{proof}

\begin{fct}$\aleph_0=\aleph_0^2=\Fin (\aleph_0)=\Fin^2(\aleph_0)= \Seq
(\aleph_0)=\seq (\aleph_0)<\fzw^{\aleph_0}$\end{fct}

\begin{proof} The only non-trivial part is $\aleph_0<\fzw^{\aleph_0}$, which
follows by the Cantor Theorem. \end{proof}

Three non-trivial relationships are given in the following

\begin{prop}\label{prop:three}
For any infinite cardinal $\fm$ we have:\hfill\\ \hspace*{3ex}(1)
$\Fin (\fm)<\fzw^\fm$\hfill\\ \hspace*{3ex}(2) $\Seq (\fm)\neq
\fzw^\fm$\hfill\\ \hspace*{3ex}(3) $\seq (\fm)\neq \fzw^\fm$
\end{prop}

(These three relationships are proved in \cite{LoriShelah}.)

\section{Permutation models}\label{sec:permutation-models}

In this section we give the definition of permutation models ({\sl
cf.}\;also~\cite{Jechchoice}). We will use permutation models to
derive relative consistency results. But first we have to introduce
models of ZFA, which is set theory with atoms ({\sl
cf.}\;\cite{Jechchoice}). Set theory with atoms is characterized by
the fact that it admits objects other than sets, namely {\bf atoms},
(also called {\bf urelements}). Atoms are objects which do not have
any elements but which are distinct from the empty-set. The
development of the theory ZFA is very much the same as that of ZF
(except for the definition of ordinals, where we have to require that
an ordinal does not have atoms among its elements). Let $S$ be a set,
then by transfinite recursion on $\alpha\in\On$ we can define
${\cP}^\alpha(S)$ as follows: ${\cP}^\emptyset(S)
:=S$, ${\cP}^{\alpha
+1}(S):={\cP}^\alpha(S)\cup{\cP}({\cP}^\alpha(S))$ and
${\cP}^\alpha(S):=\bigcup_{\beta\in\alpha}{\cP}^\alpha(S)$ when
$\alpha$ is a limit ordinal. Further let
${\cP}^\infty(S):=\bigcup_{\alpha\in\On}{\cP}^\alpha(S)$. If
${\mathcal M}$ is a model of ZFA and $A$ is the set of atoms of
${\mathcal M}$, then we have ${\mathcal M}:={\cP}^\infty(A)$. The
class $M_0:={\cP}^\infty(\emptyset)$ is a model of ZF and is called
the {\bf kernel}. Note that all the ordinals are in the kernel.

The underlying idea of permutation models, which are models of ZFA,
is the fact that the axioms of ZFA do not distinguish between the
atoms, and so a permutation of the set of atoms induces an
automorphism of the universe. The method of permutation models was
introduced by Adolf~Fraenkel and, in a precise version (with
supports), by Andrzej~Mostowski. The version with filters is due to
Ernst~Specker in \cite{Specker2}.

In the permutation models we have a set of atoms $A$ and a group
$\cG$ of permutations (or automorphisms) of $A$ (where a permutation
of $A$ is a one-to-one mapping from $A$ onto $A$). We say that a set
$\cF$ of subgroups of $\cG$ is a {\bf normal filter} on $\cG$ if for
all subgroups $H,K$ of $\cG$ we have:

\noindent \hspace*{3ex}(A) $\cG\in\cF$;\hfill\\ \hspace*{3ex}(B) if
$H\in\cF$ and $H\subs K$, then $K\in\cF$;\hfill\\ \hspace*{3ex}(C) if
$H\in\cF$ and $K\in\cF$, then $H\cap K\in\cF$;\hfill\\
\hspace*{3ex}(D) if $\pi\in\cG$ and $H\in\cF$, then $\pi
H\pi^{-1}\in\cF$; \hfill\\ \hspace*{3ex}(E) for each $a\in A$,
$\{\pi\in\cG :\pi a\,=\,a\}\in\cF$. \hfill\medskip

Let $\cF$ be a normal filter on $\cG$. We say that $x$ is {\bf
symmetric} if the group $$\Sym_{\cG}(x):=\{\pi \in\cG :\pi
x\,=\,x\}$$ belongs to $\cF$. By (E) we have that every $a\in A$ is
symmetric. \hfill\smallskip

Let $\cV$ be the class of all hereditarily symmetric objects, then
${\cV}$ is a transitive model of ZFA. We call ${\cV}$ a permutation
model. Because every $a\in A$ is symmetric, we get that the set of
atoms $A$ belongs to ${\cV}$.

Now every $\pi\in\cG$ induces an $\in$-automorphism of the universe
${\cV}$, which we denote by $\hat{\pi}$ or just $\pi$.

Because $\emptyset$ is hereditarily symmetric and for all ordinals
$\alpha$ the set $\cP^\alpha (\emptyset)$ is hereditarily symmetric
too, the class $V:={\cP}^\infty (\emptyset)$ is a class in $\cV$
which is equal to the kernel $M_0$.

\begin{fct}\label{fct:kernel}
For any ordinal $\alpha$ and any $\pi\in\cG$ we have
$\pi \alpha=\alpha$.
\end{fct}

(This one can see by induction on $\alpha$, where
$\pi \emptyset=\emptyset$ is obvious.)

Since the atoms $x \in A$ do not contain any elements, but are
distinct from the empty-set, the permutation models are models of ZF
without the axiom of foundation. However, with the Jech-Sochor
Embedding Theorem ({\sl cf.}\;\cite{JechSochor}, \cite{Jechchoice}
or~\cite{Jechbook}) one can embed arbitrarily large fragments of a
permutation model in a well-founded model of ZF:

{\sc{Jech-Sochor Embedding Theorem:}} {Let ${\mathcal M}$ be a model
of ZFA + AC, let $A$ be the set of all atoms of ${\mathcal M}$, let
${M}_0$ be the kernel of $\mathcal M$ and let $\alpha$ be an ordinal
in $\mathcal M$. For every permutation model ${\cV}\subs{\mathcal M}$
(a model of ZFA) there exists a symmetric extension $V\supseteq M_0$
(a model of ZF) and an embedding $x\mapsto\tilde{x}$ of ${\mathcal
V}$ in $V$ such that $$ (\cP^\alpha (A))^{\cV}\text{\ is
$\in$-isomorphic to\ }(\cP^\alpha (\tilde{A}))^{V}\,. $$}

Most of the well-known permutation models are of the following simple
type: Let $\cG$ be a group of permutations of $A$. A family $I$ of
subsets of $A$ is a {\bf normal ideal} if for all subsets $E,F$ of
$A$ we have:

\noindent \hspace*{3ex}(a) $\emptyset\in I$;\hfill\\ \hspace*{3ex}(b)
if $E\in I$ and $F\subs E$, then $F\in I$;\hfill\\ \hspace*{3ex}(c)
if $E\in I$ and $F\in I$, then $E\cup F\in I$;\hfill\\
\hspace*{3ex}(d) if $\pi\in\cG$ and $E\in I$, then $\pi E\in
I$;\hfill\\ \hspace*{3ex}(e) for each $a\in A$, $\{a\}\in I$.
\hfill\medskip

For each set $S\subs A$, let $$\Fix_{\cG}(S):=\{\pi\in\cG:\pi
s\,=\,s\ \mbox{for all }s\in S\};$$ and let $\cF$ be the filter on
$\cG$ generated by the subgroups $\{\Fix_{\cG}(E): E\in I\}$. Then
$\cF$ is a normal filter. Further, $x$ is symmetric if and only if
there exists a set of atoms $E_{x}\in I$ such that
$$\Fix_{\cG}(E_x)\subs\Sym_{\cG}(x).$$ We say that $E_x$ is a {\bf
support} of $x$.

\section{Consistency results derived from a few permutation models}
\label{sec:consistent}

In this section we will give some relationships between the cardinals
defined in section\;\ref{sec:cardinals} which are consistent with ZF.
We will do this by investigating the relations between certain sets
in a few permutation models. Let $\cV$ be a permutation model with
the set of atoms $A$ and let $m$ be a set in $\cV$. Let
$\mathfrak{C}(m):=\{x\in\cV: \cV\models x\lesseq m\lesseq x\}$, then
$\mathfrak{C}(m)$ is a class in $\cV$. The cardinality of $m$ in the
model $\cV$ (denoted by $\fm$) is defined by
$\fm:=\mathfrak{C}(m)\cap{\cP}^\alpha(A)\cap \cV$, where $\alpha$ is
the smallest ordinal such that
$\mathfrak{C}(m)\cap{\cP}^\alpha(A)\cap \cV\neq\emptyset$. Note that
if $m$ and $n$ are two arbitrary sets in a permutation model $\cV$
and we have for example $\cV\models m\lesseq n\not\lesseq m$ (and
therefore $\cV\models \fm<\fn$), then by the Jech-Sochor Embedding
Theorem there exists a well-founded model $V$ of ZF such that
$V\models \tilde{m}\lesseq\tilde{n}\not\lesseq\tilde{m}$ and
therefore $V\models\fm<\fn$, where $\fm$ and $\fn$ are the
cardinalities of the sets $\tilde{m}$ and $\tilde{n}$. Hence, since
every relation between sets in a permutation model can be translated
to a well-founded model, to prove that a relation between some
cardinals is consistent with ZF, it is enough to find a permutation
model in which the desired relation holds between the corresponding
sets. In the sequel we will frequently make use of this method
without always mention it.

\subsection{The basic Fraenkel model}\label{sec:Fraenkel}

First we present the basic Fraenkel model ({\sl
cf.}\;\cite{Jechchoice}).

Let $A$ be a countable infinite set (the atoms), let $\cG$ be the
group of all permutations of $A$ and let $I_{\Fin}$ be the set of all
finite subsets of $A$. Obviously, $I_{\Fin}$ is a normal ideal.

Let ${\cV}_F$ ($F$ for Fraenkel) be the corresponding permutation
model, the so called {\bf basic Fraenkel model}. Note that a set $x$
is in ${\cV}_F$ iff $x$ is symmetric and each $y\in x$ belongs to
${\cV}_F$, too.

Now we will give two basic facts involving subsets of $A$.

\begin{slm}\label{lm:support1}
Let $E\in I_{\Fin}$, then each $S\subs A$ with support $E$ is either
finite or co-finite (which means $A\setminus S$ is finite). Further,
if $S$ is finite, then $S\subs E$; and if $S$ is co-finite, then
$A\setminus S\,\subs E$.
\end{slm}

\begin{proof} Let $S\subs A$ with support $E$. Because $E$ is a support of
$S$, for all $\pi\in\Fix (E)$ and every $a\in A$ we have $\pi a\in S$
if and only if $a\in S$. If $S$ is neither finite nor co-finite, the
sets $(A\setminus E)\setminus S$ and $(A\setminus E)\cap S$ are both
infinite and hence we find a $\pi\in\Fix (E)$ such that for some
$s\in S$, $\pi s\notin S$. Now, if $S$ is finite, then $S$ must be a
subset of $E$ because otherwise we have $S\setminus E\neq\emptyset$
and we find again a $\pi\in\Fix (E)$ such that for some $s\in S$,
$\pi s\notin S$. The case when $S$ is co-finite is similar.
\end{proof}

\begin{slm}\label{lm:n-Fraenkel-transfinite}
Let $A$ be the set of atoms of the basic Fraenkel model and
let $\fm$ denote its cardinality,
then ${\cV}_F\models \aleph_0 \not\le \fzw^\fm$.
\end{slm}

\begin{proof} Assume there exists a one-to-one function $f:\N\to \cP (A)$
which belongs to ${\cV}_F$. Then, because $f$ is symmetric, there
exists a finite set $E_f\subs A$ (a support of $f$) such that
$\Fix_{\cG}(E_f)\subs \Sym_{\cG}(f)$. Now let $n\in\N$ be such that
$\Fix_{\cG}(f(n))\nsubseteq\Fix_{\cG}(E_f)$ and let
$\pi\in\Fix_{\cG}(E_f)$ be such that $\pi f(n)\neq f(n)$. With the
fact\;\ref{fct:kernel} we get that $\pi n=n$ and therefore $f(\pi
n)=f(n)$. So, $E_f$ cannot be a support of $f$, which implies that
the function $f$ does not belong to ${\cV}_F$. \end{proof}

The following proposition gives the relationships in the basic
Fraenkel model between some of the cardinals defined in
section\;\ref{sec:cardinals}, where $\fm$ denotes the cardinality of
the set of atoms of ${\cV}_F$.

\begin{sprop}\label{prop:Fraenkel-model}
Let $\fm$ denote the cardinality of the set of atoms $A$ of
${\cV}_F$. Then the in the model ${\cV}_F$ we have the
following:\hfill\\ \hspace*{3ex}(1) $\Fin (\fm)\incomp \Seq
(\fm)$\hfill\\ \hspace*{3ex}(2) $\Fin (\fm)\incomp \seq
(\fm)$\hfill\\ \hspace*{3ex}(3) $\Seq (\fm)\incomp \fzw^\fm$\hfill\\
\hspace*{3ex}(4) $\seq (\fm)\incomp \fzw^\fm$
\end{sprop}

\begin{proof} (1) Assume first that there exists a function $f\in {\cV}_F$
from $\Fin (A)$ into $\Seq (A)$ and let $E_f\in I_{\Fin}$ be a
support of $f$. Choose two arbitrary distinct elements $a_0$ and
$a_1$ of $A\setminus E_f$ such that $U:=\{x\in A: x \text{\:occurs in
$f(\{a_0,a_1\}\cup E_f)$}\}\nsubseteq E_f$ and put
$E_f^{*}:=\{a_0,a_1\}\cup E_f$. Choose a $y\in U\setminus E_f$ and a
permutation $\pi\in\Fix_{\cG}(E_f)$ such that $\pi y\neq y$ and $\pi
a_i=a_{1-i}$ (for $i\in\{0,1\}$). Now, $\pi E_f^{*}=E_f^{*}$ but $\pi
f(E_f^{*})\neq f(E_f^{*})$, which implies either that $f$ is not a
function or that $E_f$ is not a support of $f$. In both cases we get
a contradiction to our assumption.\hfill\\ The fact that $\Seq
(\fm)\nleq \Fin (\fm)$ we get by $\Fin (\fm)< \fzw^{\fm}$ (see
Proposition\;\ref{prop:three}\;(1)) and  by $\Seq (\fm)\nleq
\fzw^{\fm}$ (which will be shown in (3)).

\hspace*{1.5ex}(2) Because $\Seq (\fm )\le \seq (\fm)$, by (1) it
remains to show that $\Fin (A)\not\lesseq\seq (A)$. Assume there
exists a function $g\in\cV$ from $\Fin (A)$ into $\seq (A)$ and let
$E_g\in I_{\Fin}$ be a support of it.\hfill\\ \hspace*{1ex}{--} If
for each $p\in [A\setminus E_g ]^2$ we have $\Fix_{\cG}
(E_g)\subs\Sym_{\cG}(g(p))$, then we find $\{a_0,a_1\}$ and
$\{b_0,b_1\}$ in $[A\setminus E_g ]^2$ with
$\{a_0,a_1\}\cap\{b_0,b_1\}=\emptyset$, and a permutation
$\pi\in\Fix_{\cG}(E_g)$ such that $\pi a_i=b_i$ and $\pi b_i= a_i$
(for $i\in\{0,1\}$). Now we get $\pi g(\{a_0,a_1\})= g(\{a_0,a_1\})$
and $\pi \{a_0,a_1\}=\{b_0,b_1\}$, which contradicts our
assumption.\hfill\\ \hspace*{1ex}{--} Otherwise, there exists a set
$\{a_0,a_1\}\in [A\setminus E_g ]^2$ with $\Fix_{\cG}
(E_g)\nsubseteq\Sym_{\cG}(g(\{a_0,a_1\}))$, hence we find in the
sequence $g(\{a_0,a_1\})$ an element $y\in A$ which does not belong
to $E_g$. Now let $\pi\in\Fix_{\cG}(E_g)$ be such that $\pi
a_i=a_{1-i}$ (for $i\in\{0,1\}$) and $\pi y\neq y$, then $\pi
g(\{a_0,a_1\})\neq g(\{a_0,a_1\})$ and $\pi \{a_0,a_1\}=\{a_0,a_1\}$,
which contradicts again our assumption.

\hspace*{1.5ex}(3) Because $\fm$ is infinite we have (by
Proposition\;\ref{prop:three}\;(1)) $\Fin (\fm)<\fzw^{\fm}$, which
implies (by (1)) that $\fzw^{\fm}\nleq\Seq (\fm)$ and it remains to
show that $\Seq (A)\not\lesseq {\cP}(A)$. Assume there exists a
function $h\in\cV_F$ from $\Seq(A)$ into $\cP(A)$ and let $E_h\in
I_{\Fin}$ be a support of $h$ with $|E_h|\ge 4$. Consider $\Seq
(E_h)$, then, because $|E_h|\ge 4$, it is easy to compute that $|\Seq
(E_h)|>2\cdot 2^{|E_h|}$, which implies (by Lemma\;\ref{lm:support1})
that there exists an $s_0\in\Seq (E_h)$ such that $E_h$ is not a
support of $h(s_0)$. Let $E_0:=\bigcap\{E\in I_{\Fin}:\text{$E$ is a
support of $h(s_0)$}\}$, then $E_0$ is a support of $h(s_0)$, too.
Choose a $y\in E_0\setminus E_h$ and a permutation
$\pi\in\Fix_{\cG}(E_h)$ such that $\pi y\neq y$. Now, because
$\pi\in\Fix_{\cG}(E_h)$ and $s_0\in\Seq (E_h)$ we have $\pi s_0=
s_0$, and by construction we get $\pi h(s_0)\neq h(s_0)$. This
implies either that $h$ is not a function or that $E_h$ is not a
support of $h$ and in both cases we get a contradiction to our
assumption.

\hspace*{1.5ex}(4) By $\Fin (\fm)<\fzw^{\fm}$ and
$\Fin (\fm)\nleq \seq (\fm)$ we get $\fzw^{\fm}\nleq \seq (\fm)$, and
the inequality $\seq (\fm) \nleq \fzw^{\fm}$ follows from
$\Seq (\fm )\nleq \fzw^{\fm}$ and $\Seq (\fm)\le\seq (\fm)$.
\end{proof}

\subsection{The ordered Mostowski model}\label{sec:Mostowski}

Now we shall construct the ordered Mostowski model ({\sl
cf.}\;also~\cite{Jechchoice}).

Let the infinite set of atoms $A$ be countable, and let $<^M$ be a
linear order on $A$ such that $A$ is densely ordered and does not
have a smallest or greatest element (thus $A$ is isomorphic to the
rational numbers). Let $\cG$ be the group of all order-preserving
permutations of $A$, and let again $I_{\Fin}$ be the ideal of the
finite subsets of $A$.

Let ${\cV}_M$ ($M$ for Mostowski) be the corresponding permutation
model (given by $\cG$ and $I_{\Fin}$), the so called {\bf ordered
Mostowski model}.

Because all the sets in the ordered Mostowski model are symmetric,
each subset of $A$ has a finite support. By similar arguments as in
the proof of Lemma\;\ref{lm:n-Fraenkel-transfinite} one can show

\begin{slm}\label{lm:n-Mostowski-transfinite}
Let $A$ be the set of atoms of the ordered Mostowski model and let
$\fm$ denote its cardinality, then ${\cV}_M\models \aleph_0 \not\le
\fzw^\fm$.
\end{slm}

For a finite set $E\subs A$, one can give a complete description of
the subsets of $A$ with support $E$ and one gets the following

\begin{sfct}\label{fct:Mostowski-support}
If $E\subs A$ is a finite set of cardinality $n$, then there are
$2^{2n+1}$ sets $S\subs A$ (in ${\cV}_M$) such that $E$ is a support
of $S$.
\end{sfct}

(For a proof see~\cite[p.\,32]{LoriShelah}.)

In the following we investigate the relationships between some of the
cardinals defined in section\;\ref{sec:cardinals} in the ordered
Mostowski model, where $\fm$ will be cardinality of the set of atoms
of ${\cV}_M$.

Let $\fm$ denote the cardinality of the set of atoms $A$ (of the
ordered Mostowski model). In Theorem\;1 of \cite{LoriShelah} it is
shown that $\fzw^\fm \le^* \Fin (\fm)$. Now, by Fact\;\ref{fct:*}, we
get $\fzw^{\fzw^\fm}\le \fzw^{\Fin (\fm)}$ which implies (by the
Cantor-Bernstein Theorem, as $\Fin (\fm)\le \fzw^\fm$) that the
equation $\fzw^{\fzw^\fm}= \fzw^{\Fin (\fm)}$ holds in the ordered
Mostowski model.

Unlike in the basic Fraenkel model, all the simple cardinalities
defined in section\;\ref{sec:cardinals} are comparable in the ordered
Mostowski model:

\begin{sprop}\label{prop:Mostowski-model}
Let $\fm$ denote the cardinality of the set of atoms of ${\cV}_M$.
Then the following holds in ${\cV}_M$: $$\fm < \Fin (\fm) < \fzw^\fm
< \Seq (\fm) < \seq (\fm)\,.$$
\end{sprop}

\begin{proof} Let $A$ be the set of atoms $A$ of the ordered Mostowski
model.

$\fm <\Fin (\fm)$: It is obvious that the function $f:A\to\Fin (A)$,
defined by $f(a):=\{a\}$, is a one-to-one function from $A$ into
$\Fin (A)$. Now assume that there exists also a one-to-one function
$g$ from $\Fin (A)$ into $A$. Let $a_0:=g(\emptyset)$ and
$a_{n+1}:=g(\{a_0,\ldots,a_n\})$ (for $n\in\N$). The
$\omega$-sequence $\la a_0,a_1,\ldots,a_n,\ldots\ra$ is a one-to-one
sequence of $A$, which implies that $\aleph_0\le \fm$, but this is a
contradiction to Lemma\;\ref{lm:n-Mostowski-transfinite}.

$\Fin (\fm) < \fzw^\fm$: Because $A$ is infinite, by
Proposition\;\ref{prop:three}\;(1) we have $\Fin (\fm) < \fzw^\fm$.

$\fzw^\fm <\Seq (\fm)$: For a set $S\subs A$, let $\supp
(S):=\bigcap\{E\in I_{\Fin}: E$ is a support of $S\}$, then $\supp
(S)$ is a support of $S$, too; in fact, it is the smallest support of
$S$. Using the order-relation ``$<^M$'' on the set of atoms $A$, we
can define an ordering on the set of finite subsets of $A$ as
follows. For two finite sets $\{a_0,\ldots, a_n\}$ and
$\{b_0,\ldots,b_m\}$ of $A$, where $a_i<^M a_{i+1}$ and $b_j<^M
b_{j+1}$ (for $i<n$ and $j<m$), let
$\{a_0,\ldots,a_n\}{<}_{\Fin}\{b_0,\ldots,b_m\}$ iff either $n<m$ or
for $n=m$ we have $\exists i\le n\forall j<i (a_j=b_j\wedge a_i<^M
b_i)$). The ordering ``${<}_{\Fin}$'' on the finite subsets of $A$
induces an ordering on the power-set of $A$ (because every subset of
$A$ has a well-defined smallest finite support). Further, the
order-relation ``$<^M$'' induces in a natural way an ordering on the
set of all permutations of a given finite subset of $A$ and we
identify a permutation $\tau$ of a finite subset $\{c_0<^M\ldots
<^M c_{n-1}\}$ with $\la\tau (c_0),\tau (c_1),\ldots,\tau
(c_{n-1})\ra\in\Seq (A)$. Now we choose $20$ distinct atoms $c_0<^M
c_1<^M \ldots <^M c_{19}$ of $A$ and define a function $f$ from $\cP
(A)$ into $\Seq (A)$ as follows. For $S\subs A$ with $|\supp (S)|\ge
11$, let $f(S)$ be the $k$th permutation of $\supp (S)$, where $S$ is
the $k$th subset of $A$ with smallest support $\supp (S)$ (this we
can do because for $|\supp (S)|\ge 11$ we have $|\supp (S)|\, ! \ge
2^{2 |\supp (S)|+1}$). If $\supp (S)= \{a_0,\ldots,a_l\}$ for $l\le
9$ (where $a_i<a_{i+1}$), then we choose the first 10 elements (with
respect to $<^M$) of $\{c_0,\ldots,c_{19}\}$ which are not in $\supp
(S)$, say $\{d_0,\ldots,d_9\}$ and put $f(S)=\la a_0,
\ldots,a_l,d_{\iota_0},\ldots,d_{\iota_9}\ra$, where
$d_{\iota_0}\ldots d_{\iota_9}$ is the $(10 ! -k)$th permutation of
$d_0\ldots d_9$ and $S$ is the $k$th subset of $A$ with smallest
support $\supp (S)$. By Lemma\;\ref{fct:Mostowski-support}, the
function $f$ is a well-defined one-to-one function from $\cP (A)$
into $\Seq (A)$. If there exists a one-to-one function from $\Seq
(A)$ into $\cP (A)$, then, because $n ! >2^{2n+1}+2$ for $n\ge 10$,
we can build an one-to-one $\omega$-sequence of $A$, which is a
contradiction to Lemma\;\ref{lm:n-Mostowski-transfinite}.

$\Seq(\fm)<\seq(\fm)$: Because each one-to-one sequence of $A$ is a
sequence of $A$, we have $\Seq (A)\lesseq\seq(A)$. Now assume that
there exists also a one-to-one function $g$ from $\seq (A)$ into
$\Seq (A)$. Choose an arbitrary atom $a\in A$ and let $s_n:=g(\la
a,a,\ldots,a\ra_n)$, where $\la a,a,\ldots,a\ra_n$ denotes the
sequence of $\{a\}$ of length $n$. Because for every $n\in\N$, the
sequence $s_n$ is a one-to-one sequence of $A$, for every $n\in\N$
there exists a $k>n$ and a $b\in A$ such that $b$ occurs in $s_k$ but
for $i\le n$, $b$ does not occur in $s_i$. Because a sequence is an
ordered set, with the function $g$ we can build an one-to-one
$\omega$-sequence of $A$, which contradicts
Lemma\;\ref{lm:n-Mostowski-transfinite}. \end{proof}

Let again $\fm$ denote the cardinality of the set of atoms of the
ordered Mostowski model. Using some former facts and some
arithmetical calculations, by similar arguments as is the proof of
Proposition\;\ref{prop:Mostowski-model} one can show that the
following sequence of inequalities holds in the ordered Mostowski
model:
\begin{multline*}
\fm < [\fm]^2 < \fm^2 < \Fin (\fm) < \fzw^\fm < \Seq (\fm)
< \Fin^2 (\fm) < \Seq (\Fin(\fm)) < \\
< \Fin (\fzw^\fm)< \Fin^3 (\fm) < \Fin^4 (\fm) < \ldots < \Fin^n (\fm)
< \seq (\fm) < \fzw^{\Fin (\fm)}=\fzw^{\fzw^\fm}
\end{multline*}

\subsection{A custom-built permutation model}\label{sec:custom}

In the proof of Theorem\;2 of \cite{LoriShelah}, a permutation model
${\cV}_s$ ($s$ for sequences) is constructed in which there exists a
cardinal number $\fm$ such that ${\cV}_s\models\seq (\fm) < \Fin
(\fm)$ and hence, ${\cV}_s\models\Seq (\fm) < \Fin (\fm)$.
Specifically, $\fm$ is the cardinality of the set of atoms of
${\cV}_s$.

The set of atoms of ${\cV}_s$ is built by induction, where every atom
contains a finite sequence of atoms on a lower level. We will follow
this idea, but instead of finite sequences we will put ordered pairs
in the atoms. The model we finally get will be a model in which there
exists a cardinal $\fm$, such that $\fm^2 < [\fm]^2$ (this is in fact
a finite version of Theorem\;2 of~\cite{LoriShelah}).\hfill\smallskip

\noindent We construct by induction on $n\in\N$ the following:

\noindent \hspace*{3ex}($\alpha$) $A_0$ is an arbitrary countable
infinite set.\hfill\\ \hspace*{3ex}($\beta$) $\cG_0$ is the group of
all permutations of $A_0$.\hfill\\ \hspace*{3ex}($\gamma$)
$A_{n+1}:=A_n\dot\cup\big{\{}(n+1,p,\varepsilon): p\in A_n\times
A_n\,\wedge\,\varepsilon\in\{0,1\}\big{\}}$.\hfill\\
\hspace*{3ex}($\delta$) $\cG_{n+1}$ is the subgroup of the group of
permutations of $A_{n+1}$ containing all permutations $h$ such that
for some $g_h\in\cG_n$ and $\varepsilon_h\in\{0,1\}$ we have
\begin{equation*}
h(x)=\begin{cases}
       g_h(x) & \text{if $x\in A_n$,}\\
       (n+1,g_h(p),\varepsilon_h+_2\varepsilon_x) &
       \text{if $x=(n+1,p,\varepsilon_x)$,}
     \end{cases}
\end{equation*}
where $g_h(p)=\la g_h(p_1),g_h(p_2)\ra$ for $p=\la p_1,p_2\ra$ and
$+_2$ is the addition modulo $2$.

Let $A:=\bigcup \{A_n:n\in\N\}$ and let $\operatorname{Aut}(A)$ be
the group of all permutations of $A$; then
$$\cG:=\{H\in\operatorname{Aut}(A):\forall n\in\N (H
|_{A_n}\in\cG_n\}$$ is a group of permutations of $A$. Let $\cF$ be
the normal filter on $\cG$ generated by $\{\Fix_{\cG}(E):E\subs
A\text{\ is finite}\}$, and let ${\cV}_p$ ($p$ for pairs) be the
class of all hereditarily symmetric objects.

Now we get the following

\begin{sprop}
Let $\fm$ denote the cardinality of the set of atoms $A$ of
${\cV}_p$. Then  we have ${\cV}_p\models \fm^2 < [\fm]^2$.
\end{sprop}

\begin{proof} First we show that ${\cV}_p\models \fm^2 \leq [\fm]^2$. For
this it is sufficient to find a one-to-one function $f\in{\cV}_p$
from $A^2$ into $[A]^2$. We define such a function as follows. For
$x,y\in A$ where $x=(n,p_x,\varepsilon_x)$ and
$y=(m,p_y,\varepsilon_y)$ let $$f(\la x,y\ra):=\big{\{}(n+m+1,\la
x,y\ra,0), (n+m+1,\la x,y\ra,1)\big{\}}\,.$$ For any $\pi\in\cG$ and
$x,y\in A$ we have $\pi f(\la x,y\ra)=f(\la \pi x,\pi y\ra)$ and
therefore, the function $f$ is as desired and belongs to ${\cV}_p$.

Now assume that there exists a one-to-one function $g\in{\cV}_p$ from
$[A]^2$ into $A^2$ and let $E_g$ be a finite support of $g$. Without
loss of generality we may assume that if $(n+1,\la
x,y\ra,\varepsilon)\in E_g$, then also $x,y\in E_g$. Let $k:=|E_g|$
and for $x,y\in A$ let $g(\{x,y\})=\la
t^0_{\{x,y\}},t^1_{\{x,y\}}\ra$. Let $r:=k+4$ and let
$N:=\text{Ramsey}(2,r^2,3)$, where $\text{Ramsey}(2,r^2,3)$ is the
least natural number such that for every coloring
$\tau:[\text{Ramsey}(2,r^2,3)]^2\to r^2$ we find a $3$-element subset
$H\subs \text{Ramsey}(2,r^2,3)$ such that $\tau|_{[H]^2}$ is
constant. (If $p,r,m$ are natural numbers such that $p\le m$ and
$r>0$, then by the Ramsey Theorem ({\sl
cf.}\;\cite[Theorem\;B]{Ramsey}), $\text{Ramsey}(p,r,m)$ is
well-defined.) Choose $N$ distinct elements $x_0,\ldots,x_{N-1}\in
A_0\setminus E_g$, let $X=\{x_0,\ldots,x_{N-1}\}$ and let $c_h$
($h<k$) be an enumeration of $E_g$. We define a coloring
$\tau:[X]^2\to r\times r$ as follows. For $\{x_i,x_j\}\in [X]^2$ such
that $i<j$ let $\tau (\{x_i,x_j\})=\la \tau_0 (\{x_i,x_j\}),\tau_1
(\{x_i,x_j\})\ra$ where for $l\in\{0,1\}$ we define
\begin{equation*}
\tau_l(\{x_i,x_j\}):=\begin{cases}
                   h & \text{if $t^l_{\{x_i,x_j\}}=c_h$},\\
                   k & \text{if $t^l_{\{x_i,x_j\}}=x_i$},\\
                   k+1 & \text{if $t^l_{\{x_i,x_j\}}=x_j$},\\
                   k+2 & \text{if $t^l_{\{x_i,x_j\}}\in A_0\setminus
                   (\{x_i,x_j\}\cup E_g)$},\\
                   k+3 & \text{if $t^l_{\{x_i,x_j\}}\in A\setminus (A_0\cup
                   E_g)$}.
                   \end{cases}
\end{equation*}

By the definition of $N$ we find $3$ elements
$x_{\iota_0},x_{\iota_1},x_{\iota_2}\in X$ with
$\iota_0<\iota_1<\iota_2$ such that for $l\in\{0,1\}$, $\tau_l$ is
constant on $[\{x_{\iota_0},x_{\iota_1},x_{\iota_2}\}]^2$. So, for
$\{x_{\iota_i},x_{\iota_j}\}\in
[\{x_{\iota_0},x_{\iota_1},x_{\iota_2}\}]^2$ with $i<j$ and for
$l\in\{0,1\}$, we are at least in one of the following cases:

$\begin{array}{rl} \text{\ {\it(1)}}&
t^l_{\{x_{\iota_i},x_{\iota_j}\}}=c_{h_0} \text{\ \;and\ \;}
t^{1-l}_{\{x_{\iota_i},x_{\iota_j}\}}=c_{h_1},\\ \text{{\it(2)}}&
t^l_{\{x_{\iota_i},x_{\iota_j}\}}=c_{h}\text{\ \;and\ \;}
t^{1-l}_{\{x_{\iota_i},x_{\iota_j}\}}=x_{\iota_i},\\ \text{{\it(3)}}&
t^l_{\{x_{\iota_i},x_{\iota_j}\}}=c_{h} \text{\ \;and\ \;}
t^{1-l}_{\{x_{\iota_i},x_{\iota_j}\}}=x_{\iota_j},\\ \text{{\it(4)}}&
t^l_{\{x_{\iota_i},x_{\iota_j}\}}=x_{\iota_i}\text{\ \;and\ \;}
t^{1-l}_{\{x_{\iota_i},x_{\iota_j}\}}=x_{\iota_j},\\ \text{{\it(5)}}&
t^l_{\{x_{\iota_i},x_{\iota_j}\}}\in A_0\setminus
(E_g\cup\{x_{\iota_i},x_{\iota_j}\}),\\ \text{{\it(6)}}&
t^l_{\{x_{\iota_i},x_{\iota_j}\}}\in A\setminus (E_g\cup A_0).
\end{array}$

If we are in case (1) or (2), then $g(\{x_{\iota_0},x_{\iota_1}\})=
g(\{x_{\iota_0},x_{\iota_2}\})$, and therefore $g$ is not a
one-to-one function. If we are in case (3), then $g$ is also not a
one-to-one function because $g(\{x_{\iota_0},x_{\iota_2}\})=
g(\{x_{\iota_1},x_{\iota_2}\})$.

If we are in case (4), let $\pi\in\Fix(E_g)$ be such that $\pi
x_{\iota_0}=x_{\iota_1}$ and $\pi x_{\iota_1}=x_{\iota_0}$. Assume
$g(\{x_{\iota_0},x_{\iota_1}\})=\la x_{\iota_0},x_{\iota_1}\ra$ (the
case when $g(\{x_{\iota_0},x_{\iota_1}\})=\la
x_{\iota_1},x_{\iota_0}\ra$ is symmetric). Then we have $\pi
\{x_{\iota_0},x_{\iota_1}\}=\{x_{\iota_0},x_{\iota_1}\}$, but $\pi
g(\{x_{\iota_0},x_{\iota_1}\})=\la x_{\iota_1},x_{\iota_2}\ra\neq \la
x_{\iota_0},x_{\iota_1}\ra$, and therefore $g$ is not a function in
${\cV}_p$.

If we are in case (5), let $l\in\{0,1\}$ be such that
$t^l_{\{x_{\iota_0},x_{\iota_1}\}}\in A_0\setminus
(E_g\cup\{x_{\iota_0},x_{\iota_1}\})$ and let
$a:=t^l_{\{x_{\iota_0},x_{\iota_1}\}}$. Take an arbitrary $a'\in
A_0\setminus (E_g\cup\{a,x_{\iota_0},x_{\iota_1}\})$ and let
$\pi\in\Fix(E_g\cup\{x_{\iota_0},x_{\iota_1}\})$ be such that $\pi
a=a'$ and $\pi a'=a$. Then we get $\pi
\{x_{\iota_0},x_{\iota_1}\}=\{x_{\iota_0},x_{\iota_1}\}$ but $\pi
g(\{x_{\iota_0},x_{\iota_1}\}) \neq g(\{x_{\iota_0},x_{\iota_1}\})$,
and therefore $g$ is not a function in ${\cV}_p$.

If we are in case (6), let $l\in\{0,1\}$ be such that
$t^l_{\{x_{\iota_0},x_{\iota_1}\}}\in A\setminus (E_g\cup A_0)$, thus
$t^l_{\{x_{\iota_0},x_{\iota_1}\}}=(n+1,p,\varepsilon)$ for some
$(n+1,p,\varepsilon)\in A$. Let
$\pi\in\Fix(E_g\cup\{x_{\iota_0},x_{\iota_1}\})$ be such that $\pi
(n+1,p,\varepsilon)=(n+1,p,1-\varepsilon)$. Then we have $\pi
\{x_{\iota_0},x_{\iota_1}\}=\{x_{\iota_0},x_{\iota_1}\}$ but $\pi
g(\{x_{\iota_0},x_{\iota_1}\}) \neq g(\{x_{\iota_0},x_{\iota_1}\})$,
and therefore $g$ is not a function in ${\cV}_p$.

So, in all the cases, $g$ is either not a function or it is not
one-to-one, which contradicts our assumption and completes the proof.
\end{proof}

\subsection{On sequences and the power-set}\label{sec:seq}

The Theorem\;2 of \cite{LoriShelah} states that the relation $\seq
(\fm)<\Fin (\fm)$ is consistent with ZF. If we consider the
permutation model ${\cV}_s$ ($s$ for sequences) constructed in the
proof of this theorem, we see that even more is consistent with ZF,
namely

\begin{sprop}\label{prop:Shelah-model}
It is consistent with ZF that there exists a cardinal number $\fm$,
such that $\Seq(\fm)< \seq(\fm)< \Fin(\fm)< \fzw^\fm$.
\end{sprop}

\begin{proof} Let $\fm$ denote the cardinality of the set of atoms of the
permutation model ${\cV}_s$ constructed in the proof of Theorem\;2 of
\cite{LoriShelah}. Then in ${\cV}_s$ we have $\Seq(\fm)< \seq(\fm)<
\Fin(\fm)< \fzw^\fm$:

The inequality $\seq(\fm)< \Fin(\fm)$ is Theorem\;2 of
\cite{LoriShelah} and because $\fm$ is infinite, by
Proposition\;\ref{prop:three}\;(1), we also get $\Fin(\fm)<\fzw^\fm$.

To see that also $\Seq(\fm)<\seq(\fm)$ holds in ${\cV}_s$, assume
that there exists (in ${\cV}$) a one-to-one function from $\seq
(\fm)$ into $\Seq (\fm)$. Such a one-to-one function would generate a
function $f\in {\cV}$ from $\aleph_0$ into $\fm$, but because
$f$---as an element of ${\cV}$---has a finite support, this is
impossible. \end{proof}

In the remainder of this section we show that it is consistent with
ZF that there exists a cardinal number $\fm$ such that $\Seq (\fm)<
\fzw^\fm <\seq (\fm)$. For this we construct a permutation model
${\cV}_c$ ($c$ for categorical) where $\fm$ will be the cardinality
of the set of atoms of ${\cV}_c$.

Let $L$ be the signature containing the binary relation symbol
``$<$'' and for each $n\in\N$ an $n+1$-ary relation symbol $R_n$. Let
$T_0$ be the following theory:

\noindent \hspace*{3ex}($\alpha$) $<$ is a linear order,\\
\hspace*{3ex}($\beta$) for each $n\in\N$\,:
$R_n(z_0,\ldots,z_{n})\rightarrow\bigwedge\limits_{l\neq m}(z_l\neq
z_m)\,.$

Let $K=\{N:\text{$N$ is a finitely generated structure of $T_0$}\}$,
then $K\neq\emptyset$ and further we have the following fact ({\sl
cf.}\;also~\cite[p.\,325]{Hodges}).

\begin{sfct}\label{fct:amalgamation}
$K$ has the amalgamation property. \end{sfct}

\begin{proof} If $N_0\subs N_1\in K$, $N_0\subs N_2\in K$ and $N_1\cap
N_2=N_0$, then we can define an $N\in K$ such that
$\mdom(N)=\mdom(N_1)\cup\mdom(N_2)$, $N_1\subs N$, $N_2\subs N$,
$<^{N_1}\cup <^{N_2}\subs <^N$ and for any $n\in\N$ we have
$R_n^{N_1}\cup R_n^{N_2}=R_n^{N}$.\end{proof}

As a consequence of Fact\;\ref{fct:amalgamation} we get the

\begin{slm}\label{lm:age}
There exists (up to isomorphism) a unique structure $M$ of $T_0$ such
that the cardinality of $\mdom(M)$ is $\aleph_0$, each structure
$N\in K$ can be embedded in $M$ and every isomorphism between
finitely generated substructures of $M$ (between two structures of
$K$) extends to an automorphism of $M$.
\end{slm}

\begin{proof} For a proof see {\sl e.g.}\ Theorem\;7.1.2 of \cite{Hodges}.
\end{proof}

Therefore, Th($M$) is $\aleph_0$-categorical and, because every
isomorphism between finitely generated substructures of $M$ extends
to an automorphism of $M$, the structure $M$ has non-trivial
automorphisms.

Now we construct the permutation model ${\cV}_c$ as follows. The set
$\mdom(M)$ constitutes the set of atoms $A$ of ${\cV}_c$ and $\cG$ is
the group of all permutations $\pi$ of $A$ such that: $M\models x<^M
y$ iff $M\models \pi x<^M \pi y$ and for each $n\in\N$, $M\models
R_n(z_0,\ldots,z_{n})$ iff $M\models R_n(\pi z_0,\ldots,\pi z_{n})$.
In fact, the group $\cG$ is the group of all automorphisms of $M$.
Further, let $\cF$ be the normal filter on $\cG$ generated by
$\{\Fix(E)_{\cG}:\text{$E\subs A$ is finite}\}$ and let ${\cV}_c$ be
the class of all hereditarily symmetric objects.

\notation If $i+k=n$ and $\bar y=\la y_0,\ldots,y_{n-1}\ra$, then we
write $R_{i,k}(x,\bar y)$ instead of
$R_n(y_0\ldots,y_{i-1},x,y_{i},\ldots y_{i+k-1})$. If $i=0$, we write
just $R_n(x,\bar y)$.

The following lemma follows from the fact that every isomorphism
between two structures of $K$ extends to an automorphism of $M$ and
from the fact that Th($M$) is $\aleph_0$-categorical ({\sl
cf.}\;also~\cite[Theorem\;7.3.1]{Hodges}).

\begin{slm}\label{lm:c-support}
For every set $S\subs A$ in ${\cV}_c$ there exists a unique smallest
support $\supp(S)$ and for each finite set $E\subs A$, the set
$\{S\subs A:S\in{\cV}_c\wedge \supp(S)=E\}$ is finite.
\end{slm}

\begin{proof} For $n\in\N$ let $E=\{e_0,\ldots,e_{n-1}\}\subs A$ be a finite
set of atoms. Further, let $\Theta_E$ the set of all atomic
$L$-formulas $\varphi_i(x)$ such that we have $\varphi_i(x)$ is {\it
either\/} the formula $x=e_j$ (for some $j<n$) {\it or\/} $x<^M e_j$
(for some $j<n$) {\it or\/} $R_{i,k}(x,\bar e)$ (for $i+k\le n$ and
$\bar e\in\Seq(E)$). For an atom $a\in A$ let
$$\vartheta_E(a):=\{\varphi_i(x)\in\Theta_E:M\models\varphi_i(a)\}\,;$$
thus, $\vartheta_E(a)$ is the set of all atomic formulas in
$\Theta_E$ such that $\varphi_i(a)$ holds in $M$.

Take an arbitrary $S\subs A$ in ${\cV}_c$ and let $E$ be a support of
$S$. If $s,t\in A$ are such that $\vartheta_E(t)=\vartheta_E(s)$,
then we find (by construction of $M$ and $\cG$) a permutation
$\pi\in\Fix_{\cG}(E)$ such that $\pi s=t$ and therefore we have $s\in
S$ if and only if $t\in S$. Hence, the set $S$ is determined by
$\{\vartheta_E(s):s\in S\}$, which is a finite set of finite sets of
atomic formulas.

Now we show that if $E_1$ and $E_2$ are two distinct supports of a
set $S\subs A$, then $E_1\cap E_2$ is also a support of $S$. If
$E_1\subs E_2$ or $E_2\subs E_1$, then it is obvious that $E_1\cap
E_2$ is a support of $S$. So, assume that $E_1\setminus E_2$ and
$E_2\setminus E_1$ are both non-empty and let $E_0:=E_1\cap E_2$.
Take an arbitrary $s_0\in S$ and let
$\vartheta_0:=\vartheta_{E_0}(s_0)$. Let $t\in A$ be any atom such
that $\vartheta_{E_0}(t)=\vartheta_0$. We have to show that also
$t\in S$. If $x=e_j$ belongs to $\vartheta_0$ (and thus $e_j\in
E_0$), then $s_0=t=e_j$ and we have $t\in S$. So, assume that $x=e_j$
does not belong to $\vartheta_0$ (for any $e_j\in E_0$). If $x=e_i$
does not belong to $\vartheta_{E_1}(t)$, let $t':=t$. Otherwise, if
$x=e_i$ belongs to $\vartheta_{E_1}(t)$, because
$\vartheta_{E_0}(t)=\vartheta_0$ and $E_0=E_1\cap E_2$ we have
$e_i\in E_1\setminus E_2$. By construction of $M$ we find a $t'\in A$
such that $t'\notin E_1$ and
$\vartheta_{E_2}(t')=\vartheta_{E_2}(t)$, hence, $t'\in
S\Leftrightarrow t\in S$. Now let $$\Xi_2:=\{\vartheta_{E_2}(s):s\in
S\wedge\vartheta_{E_2}(s)\cap \vartheta_{E_0}(s)=\vartheta_0\}\,.$$
Because $t'\notin E_1$ we find (again by construction of $M$) a
$t''\in A$ such that $\vartheta_{E_1}(t'')=\vartheta_{E_1}(t')$ and
$\vartheta_{E_2}(t'')\in\Xi_2$. Now, by
$\vartheta_{E_2}(t'')\in\Xi_2$ we have $t''\in S$, by
$\vartheta_{E_1}(t'')=\vartheta_{E_1}(t')$ we have $t''\in
S\Leftrightarrow t'\in S$, and because $t'\in S\Leftrightarrow t\in
S$ we finally get $t\in S$.

Hence, $\supp(S):=\bigcap\{E:\text{$E$ is a support of $S$}\}$ is a
support of $S$ and by construction it is unique. \end{proof}

Now we are ready to prove the

\begin{sprop}\label{prop:berlin}
Let $\fm$ denote the cardinality of the set of atoms of ${\cV}_c$,
then we have ${\cV}_c\models\Seq(\fm)<\fzw^\fm <\seq(\fm)$.
\end{sprop}

\begin{proof} First we show ${\cV}_c\models\Seq(\fm)<\fzw^\fm$. For $\bar
y=\la y_0,\ldots,y_{n-1}\ra\in\Seq(A)$ let $$\Phi(\bar y):\{x\in
A:M\models R_n(x,\bar y)\}\,.$$ By the construction of ${\cV}_c$, the
function $\Phi$ belongs to ${\cV}_c$ and is a one-to-one mapping from
$\Seq(A)$ to ${\cP}(A)$. Hence, ${\cV}_c\models\Seq(\fm)\le\fzw^\fm$
and because (by Proposition\;\ref{prop:three}\;(2))
$\Seq(\fm)\neq\fzw^\fm$ is provable in ZF, we get
${\cV}_c\models\Seq(\fm)<\fzw^\fm$.

To see that ${\cV}_c\models\fzw^\fm<\seq(\fm)$, notice first that by
Proposition\;\ref{prop:three}\;(3), the inequality
$\seq(\fm)\neq\fzw^\fm$ is provable in ZF, and therefore it is enough
to find a one-to-one function from ${\cP}(A)$ into $\seq(A)$ which
lies in ${\cV}_c$. For each finite set $E\subs A$, let $\eta_E$ be an
enumeration of the set $\Theta_E$. (The function $E\mapsto\eta_E$
exists as $<^M$ is a linear order on the finite set $E$.) Then by the
Lemma\;\ref{lm:c-support} and its proof, for each finite set $E\subs
A$, $\eta_E$ induces a mapping from $\{S\subs A:\supp(S)=E\}$ into
$k$, for some $k\in\N$. Now fix two distinct atoms $a,b\in A$ and let
$$\begin{array}{rccc}{\Psi}
:&{\cP}(A)&\longrightarrow &\seq(A)\\ & S& \longmapsto & \la
e_0,\ldots,e_{n-1},a,b,\ldots,b\ra\end{array}$$ be defined as
follows: $E=\{e_0,\ldots,e_{n-1}\}:=\supp(S)$ such that $e_0<^M
\ldots<^M e_{n-1}$ and the length of the sequence $\Psi(S)$ is equal
to $n+1+l$, where $\eta_E$ maps $S$ to $l$. The function $\Psi$ is as
desired, because it is a one-to-one function from ${\cP}(A)$ into
$\seq(A)$ which lies in ${\cV}_c$. \end{proof}

\rmk Because the relation $<^M$ is a dense linear order on the set of
atoms of ${\cV}_c$, with similar arguments as in the proof of the
Proposition\;\ref{prop:Mostowski-model} one can show that
${\cV}_c\models\Fin(\fm)<\Seq(\fm)$ (where $\fm$ denotes the
cardinality of the atoms of ${\cV}_c$).

\section{Cardinals related to the power-set}\label{sec:misc}

In this section we compare the cardinalities of some sets which are
related to the power-set. First we consider the power-set itself and
afterwards we give some results involving the set of partitions.

The following fact can be found also in \cite{LindenbaumTarski} or
\cite[VIII\,2\;Ex.\,9]{Sierpinski}. However, we want to give here a
combinatorial proof of this fact.

\begin{fct}\label{fct:aleph0}
If $\aleph_0\le \fzw^\fm$, then $\fzw^{\aleph_0}\le \fzw^\fm$.
\end{fct}

\begin{proof} Take an arbitrary $m\in\fm$. Because $\aleph_0\le\fzw^\fm$, we
find an one-to-one $\omega$-sequence $\la
p_0,p_1,\ldots,p_n,\ldots\ra$ of ${\cP}(m)$. Define an equivalence
relation on $m$ by $$x\sim y\ \text{if and only if}\ \forall n\in\N
(x\in p_n\leftrightarrow y\in p_n)\,,$$ and let $[x]:=\{y\in m:y\sim
x\}$. For $x\in m$ let $g[x]:=\{n\in\N: x\in p_n\}$, then, for every
$x\in m$, we have $g[x]\subs\N$ and $g[x]=g[y]$ if and only if $[x]=
[y]$. We can consider $g[x]$ as an $\omega$-sequence of $\{0,1\}$ by
stipulating $g[x](n)=0$ if $x\in p_n$ and $g[x](n)=1$ if $x\notin
p_n$. Now we define an ordering on the set $\{g[x]:x\in m\}$ as
follows: \begin{multline*}g[x]<_g g[y]\hspace*{2ex}\text{\it{if and
only if}}\\ \exists n\in\N\big(g[x](n)<g[y](n)\wedge\forall
k<n(g[x](k)=g[y](k))\big)\,.\end{multline*} This is a total order on
the set $\{g[x]:x\in m\}$. Let $P_n^0:=\{g[x]:g[x](n)=0\}$, then for
each $n\in\N$ the set $P_n^0$ is a set of $\omega$-sequences of
$\{0,1\}$. The order relation $<_g$ defines an ordering on each
$P_n^0$ and we must have one of the following two cases:

\noindent \hspace*{2ex}{\it Case\;1\/}: For each $n\in\N$, $P_n^0$ is
well-ordered by $<_g$. \hfill\\ \hspace*{2ex}{\it Case\;2\/}: There
exists a least $n\in\N$ such that $P_n^0$ is not well-ordered by the
relation $<_g$.

If we are in case\;1, then we find a well-ordering on
$\bigcup_{n\in\N} P_n^0$. Let the ordinal $\alpha$ denote its
order-type, then $\alpha\ge\omega$ (otherwise the $\omega$-sequence
$\la p_0,p_1,\ldots\ra$ would not be one-to-one) and therefore we can
build a one-to-one $\omega$-sequence $\la g[x_0],g[x_1],\ldots\ra$ of
$\{g[x]:x\in m\}$. If we define $q_i:=\{x\in m:g[x]=g[x_i]\}$, then
the set $Q:=\{q_i:i\in\N\}$ is a set of pairwise disjoint subsets of
$m$ of cardinality $\aleph_0$. Therefore, the cardinality of
${\cP}(Q)$ is $\fzw^{\aleph_0}$ and because for $q\subs Q$ the
function $\varphi(q):=\bigcup q\subs m$ is a one-to-one function, we
get $\fzw^{\aleph_0}\le\fzw^\fm$.

If we are in case\;2, let $n$ be the least natural number such that
$P_n^0$ is not well-ordered by $<_g$. Let $S_0:=\bigcup\{s\subs
P_n^0: s\;\text{has no smallest element}\}$. Then $S_0\subs P_n^0$
has no smallest element, too. For $k\in\N$ we define $S_{k+1}$ as
follows. If $S_k\cap P_{n+k+1}^0=\emptyset$, then $S_{k+1}:= S_k$;
otherwise, $S_{k+1}:=S_k\cap P_{n+k+1}^0$. By construction, for every
$k\in\N$, the set $S_k$ is not the empty set and it is not
well-ordered by $<_g$. Thus, for every $k\in\N$ there exists an $l>k$
such that $S_l$ is a proper subset of $S_k$. Now let $\la
S_{k_0},S_{k_1},\ldots\ra$ be such that for all $i<j$ we have
$S_{k_i}\setminus S_{k_j}\neq\emptyset$ and let $q_i:=\{x\in m:
g[x]\in (S_{k_i}\setminus S_{k_{i+1}}\}$. Then the set
$Q:=\{q_i:i\in\N\}$ is again a set of pairwise disjoint subsets of
$m$ of cardinality $\aleph_0$ and we can proceed as above.
\end{proof}

\begin{fct}\label{fct:n-mal}
If $\aleph_0\nleq\fzw^\fm$, then for every natural number $n$ we have
$n\cdot\fzw^\fm<(n+1)\cdot\fzw^\fm$ and if $\emptyset\notin\fm$ we
also have $\fzw^{n\cdot\fm}<\fzw^{(n+1)\cdot\fm}$.
\end{fct}

\begin{proof} We will give the proof only for the former case, since the
proof of the latter case is similar. Let $n$ be an arbitrary natural
number. It is obvious that we have $n\cdot\fzw^\fm\le
(n+1)\cdot\fzw^\fm$. So, for an $m\in\fm$, let us assume that we also
have a one-to-one function $f$ from $(n+1)\times{\cP}(m)$ into
$n\times{\cP}(m)$. For $k\ge 1$ let $\la s_0,\ldots,s_{k-1}\ra_k$ be
a one-to-one $k$-sequence of $\cP(m)$ and let $U_k:=\{s_i: i<k\}$. We
can order the set $(n+1)\times U_k$ as follows: $\la l_i,s_i\ra<_U
\la l_j,s_j\ra$ iff either $i<j$ or $i=j\rightarrow l_i<l_j$. Because
$\mid (n+1)\times U_k\mid = (n+1)\cdot k$ and $k\ge 1$, we have
$(n+1)\cdot k>n\cdot k$ and hence there exists a first $\la
l_i,s_i\ra$ (w.r.t.\;$<_U$), such that the second component of $f(\la
l_i,s_i\ra )$ does not belong to $U_k$. Now we define $s_k:= f(\la
l_i,s_i\ra )$ and the $(k+1)$-sequence $\la s_0,\ldots,s_k \ra_{k+1}$
is a one-to-one sequence of $\cP (m)$. Repeating this construction,
we finally get an one-to-one $\omega$-sequence of $\cP (m)$. But this
is a contradiction to $\aleph_0\not\leq\fzw^\fm$. So, our assumption
was wrong and we must have $n\cdot\fzw^\fm < (n+1)\cdot\fzw^\fm$.
\end{proof}

Because it is consistent with ZF that there exists an infinite
cardinal number $\fm$ such that $\aleph_0\nleq\fzw^\fm$ (see
Lemma\;\ref{lm:n-Fraenkel-transfinite}), it is also consistent with
ZF that there exists an infinite cardinal $\fm$ such that
$\fzw^\fm<\fzw^\fm+\fzw^\fm$. Concerning $\fzw^{\fzw^\fm}$,
Hans~L\"auchli proved in ZF that for every infinite cardinal number
$\fm$ we have $\fzw^{\fzw^\fm}+\fzw^{\fzw^\fm}=\fzw^{\fzw^\fm}$
(see~\cite{Lauchli1}). In particular, he got this result as a
corollary of the following: It is provable in ZF that for any
infinite cardinal $\fm$ we have $\big(
\fzw^{\Fin(\fm)}\big)^{\aleph_0} =\fzw^{\Fin(\fm)}$ ({\sl
cf.}\;\cite{Lauchli1}). Now, because $\fzw^\fm =\Fin(\fm)+\fq$ (for
some $\fq$), we have
$\fzw^{\fzw^\fm}=\fzw^{\Fin(\fm)+\fq}=\fzw^{\Fin(\fm)}\cdot\fzw^\fq$,
and therefore, $\fzw^{\fzw^\fm}=\big(
\fzw^{\Fin(\fm)}\big)^{\aleph_0}\cdot\fzw^\fq\ge\fzw^{\Fin(\fm)}\cdot
\fzw^{\Fin(\fm)}\cdot\fzw^\fq\ge\fzw\cdot
\fzw^{\Fin(\fm)+\fq}=\fzw^{\fzw^\fm}+\fzw^{\fzw^\fm}$, and the
equation $\fzw^{\fzw^\fm}+\fzw^{\fzw^\fm}=\fzw^{\fzw^\fm}$ follows by
the Cantor-Bernstein Theorem.

Now we give some results concerning the set of partitions of a given
set.

A set $p\subs \cP(m)$ is a partition of $m$ if $p$ is a set of
pairwise disjoint, non-empty sets such that $\bigcup p = m$. We
denote the set of all partitions of a set $m$ by $\Part (m)$ and the
cardinality of $\Part (m)$ by $\Part (\fm)$. Because each partition
of $m$ is a subset of the power set of $m$, we obviously have $\Part
(\fm)\leq \fzw^{\fzw^{\fm}}.$ It is also easy to see that if $m$ has
more than $4$ elements, then $\fzw^{\fm}\leq \Part (\fm).$ If we
assume the axiom of choice, then for every infinite set $m$ we have
$\fzw^{\fm}=\Part (\fm)$ ({\sl cf.}\;\cite[XVII.4
Ex.\,3]{Sierpinski}). But it is consistent with ZF that there exists
an infinite set $m$ such that $\fzw^{\fm}<\Part (\fm)$. Moreover we
have the following

\begin{prop}\label{fct:lower-bound}
If $\fm\ge 5$ and $\aleph_0\nleq \fzw^\fm$, then $\fzw^{\fm}< \Part
(\fm)$.
\end{prop}

\begin{proof} For a finite $\fm\ge 5$, it is easy to compute that $\fzw^\fm<
\Part(\fm)$. So, let us assume that $\fm$ is infinite and take
$m\in\fm$. Because $m$ has more than $4$ elements we have ${\cP}(m)
\lesseq \Part(m)$. Now let us further assume that there exists a
one-to-one function $f$ from $\Part (m)$ into ${\cP}(m)$. First we
choose $4$ distinct elements $a_0,a_1,a_2,a_3$ from $m$. Let $c_i:=
\{a_i\}$ (for $i<4$) and $c_4:=m\setminus\bigcup\{c_i:i<4\}$, then
$C_5:=\{c_i:i\le 4\}$ is a set of pairwise disjoint, non-empty
subsets of $m$ such that $\bigcup C_5=m$. Let $S_k=\la
X_0,\ldots,X_{k-1}\ra_k$ be a one-to-one sequence of ${\cP}(m)$ of
length $k$. With respect to the sequence $S_k$ we define an
equivalence relation on $m$ as follows. $x\sim y$ if and only if for
all $i<k$: $x\in X_i\Leftrightarrow y\in X_i$. For $x\in m$ let
$[x]:=\{y\in m:y\sim x\}$ and let $\chi_x: k\to \{0,1\}$ be such that
$\chi_x(i)=0$ if and only if $x\in X_i$. Notice that we have
$\chi_x=\chi_y$ if and only if $x\sim y$. We define an ordering on
the set of equivalence classes by stipulating $[x]<_{\chi}[y]$ if
there exists an $i<k$ such that $\chi_x(i)<\chi_y(i)$ and for all
$j<i$ we have $\chi_x(j)=\chi_y(j)$. Further let
$C_k=\cC(S_k):=\{[x]:x\in m\}$, then $C_k$ is a set of pairwise
disjoint, non-empty subsets of $m$ such that $\bigcup C_k=m$ and
$|C_k|\ge k$.

Let us assume that we already have constructed a set
$C_k=\{c_0,\ldots,\linebreak[2]c_{k-1}\}$ (for some $k\ge 5$) where
$C_k=\cC(S_k)$ and $S_k=\la X_0,\ldots,X_{k-1}\ra_k\in\Seq({\cP}(m))$
for some $k\ge 5$. Every partition of $l=|\cC_k|$ induces a partition
of $m$ (this is because of the properties of $C_k$) and hence we get
a one-to-one mapping $\iota$ from $\Part (l)$ into $\Part (m)$.
Notice that the ordering $<_{\chi}$ on $C_k$ induces an ordering on
$\Part (l)$. Because $l\ge k\ge 5$ we have $|\Part (l)|>|{\cP}(l)|$
and therefore we find a first partition $q$ of $l$ (first in the
sense of the ordering on $\Part (l)$) such that the set $f(\iota
(q))$ is not the union of elements of $C_k$. We define $X_k:=f(\iota
(q))$, $S_{k+1}:=\la X_0,\ldots,X_k\ra_{k+1}$ and
$C_{k+1}:=\cC(S_{k+1})$. Repeating this construction, we finally get
an one-to-one $\omega$-sequence of ${\cP}(m)$. But this is a
contradiction to $\aleph_0\nleq\fzw^\fm$ and therefore we have $\Part
(m)\not\lesseq {\cP}(m)$ and by ${\cP}(m)\lesseq\Part(m)$ we get
$\fzw^\fm
<\Part(\fm)$.
\end{proof}

One can consider a partition of a set $m$ also as a subset of
$[m]^2$. To see this, let $f:\Part (m)\to\cP([m]^2)$ be such that for
$p\in\Part (m)$ we have $\{i,j\}\in f(p)$ if and only if $\exists
b\in p(\{i,j\}\subs b)$. Therefore, for any cardinal $\fm$ we have
$\Part (\fm)\le \fzw^{[\fm ]^2}$ and as a consequence we get

\begin{fct}\label{fct:pair-power}
If $\fm\ge 4$ and $\aleph_0\nleq\fzw^\fm$, then
$\fzw^{\fm}<\fzw^{[\fm ]^2}$.
\end{fct}

\begin{proof} This follows from the Fact\;\ref{fct:lower-bound} and the fact
that $\Part (\fm)\le \fzw^{[\fm ]^2}$. \end{proof}

Let CH$(\fm)$ be the following statement: If $\fn$ is a cardinal
number such that $\fm\le\fn\le \fzw^\fm$, then $\fn = \fm$ or $\fn =
\fzw^\fm$. Specker showed in \cite{Specker1} that if CH$(\fm)$ holds
for every infinite cardinal $\fm$, then we have the axiom of choice.

Concerning the set of partitions we get the following easy
\begin{fct}\label{fct:strong}
If $\fm$ is infinite and ${\operatorname{CH}}(\fm)$ holds, then
$\Part (\fm)=\fzw^\fm$.
\end{fct}

\begin{proof} Note that $\fm\le [\fm]^2\le\Fin(\fm)<\fzw^\fm$ and therefore,
by ${\operatorname{CH}}(\fm)$, we must have $\fm = [\fm]^2$, and by
$\fzw^\fm\le\Part (\fm)\le \fzw^{[\fm]^2}$ we get $\Part (\fm)=
\fzw^\fm$. \end{proof}

The assumption in Fact\;\ref{fct:strong} is of course very strong.
For example it is also consistent with ZF that there exists an
infinite set $m$ such that $[\fm]^2 >\fm$ ({\sl e.g.}, let $m$ be the
set of atoms in the basic Fraenkel model or in the model ${\cV}_p$).
Moreover, in the second Cohen model constructed in
\cite[5.4]{Jechchoice}---which is a symmetric model---there exists a
set $m$ such that $\aleph_0\le [\fm]^2$, $\aleph_0\not\le\fm$ and
$\aleph_0\le^*\fm$.

One cannot expect that the cardinality of a partition $p\in\Part (m)$
is very large: If $\fp$ is a partition of $\fm$, then $\fp\le^*\fm$
and (by Fact\;\ref{fct:*}) we get $\fzw^\fp \le \fzw^\fm$, which
implies $\fp < \fzw^\fm$. On the other hand, for $p\in\Part (m)$ we
can have ${\mathfrak{p}}>\fm$. To see this, take any two cardinal
numbers $\fn$ and $\fm$ such that $\fn<\fm$ and $\fm\le^* \fn$
(examples for such cardinals can be found {\sl e.g.}\;in
\cite{LoriShelah}). Now take $m\in\fm$ and $n\in\fn$, then by the
definition of $\le^*$ there exists a function $f$ from $n$ onto $m$
and the set $p:=\big{\{}\{x\in n:f(x)=y\}:y\in m\big{\}}$ is a
partition of $n$ of cardinality $\fm$. Moreover, this can also happen
even if we partition the real line:

\begin{fct}\label{fct:partition}
It is consistent with ZF that the real line can be partitioned into a
family $p$, such that ${\mathfrak{p}}>\fzw^{\aleph_0}$, where
$\fzw^{\aleph_0}$ is the cardinality of the set of the real numbers.
\end{fct}

\begin{proof} Specker showed in \cite[II 3.32]{Specker2} that if the
real numbers are the countable union of countable sets, then
$\aleph_1$ and $\fzw^{\aleph_0}$ are incomparable. Furthermore,
Henri~Lebesgue gave in \cite{Lebesgue} a proof that
$\aleph_1\le^*\fzw^{\aleph_0}$ (see also \cite[XV\;2]{Sierpinski}).
Therefore we can decompose effectively the interval $(0,1)$ into
$\aleph_1$ disjoint non-empty sets and obtain a decomposition of the
real line into $\aleph_1 +2^{\aleph_0}$ disjoint non-empty sets. If
$\aleph_1\nleq 2^{\aleph_0}$, then $2^{\aleph_0}<\aleph_1
+2^{\aleph_0}$. Hence, in the model of Solomon~Feferman and
Azriel~Levy ({\sl cf.}\;\cite{FefermanLevy})---in which the real
numbers are the countable union of countable sets---we find a
decomposition of the real line into more than $2^{\aleph_0}$ disjoint
non-empty sets (see also \cite[p.\,372]{Sierpinski}). \end{proof}

\section{Summary}

First we summarize the results we got in the
sections\;\ref{sec:provable} and~\ref{sec:consistent} by listing all
the possible relationships between the cardinal numbers ${\fm}$,
${\Fin (\fm)}$, $\Seq (\fm)$, $\seq (\fm)$ and $\fzw^\fm$, where the
cardinal number $\fm$ is infinite.\pagebreak

{}\vspace*{16pt}

\begin{tabular}{|c|c|c|c|c|c|}
\vspace*{-16pt} {\hspace*{1.3cm}} & {\hspace*{1.3cm}} &
{\hspace*{1.3cm}} & {\hspace*{1.3cm}} & {\hspace*{1.3cm}} &
{\hspace*{1.3cm}} \\ \hline {\ } &
$\underset{\rule[-0.8ex]{0pt}{0pt}}{\fm}$ &
$\overset{\rule[2ex]{0pt}{0pt}}{\Fin (\fm)}$ & $\Seq (\fm)$ & $\seq
(\fm)$ & $\fzw^\fm$ \\ \hline $\fm$ &
$\underset{\rule[-0.5ex]{0pt}{0pt}}{\overset{\rule[2ex]{0pt}{0pt}}{=}}$
& $\overset{\scriptscriptstyle{\ref{sec:provable}}}{=}\
\overset{\scriptscriptstyle{\ref{sec:Mostowski}}}{<}$ &
$\overset{\scriptscriptstyle{\ref{sec:provable}}}{=}\
\overset{\scriptscriptstyle{\ref{sec:Mostowski}}}{<}$ &
$\overset{\scriptscriptstyle{\ref{sec:provable}}}{=}\
\overset{\scriptscriptstyle{\ref{sec:Mostowski}}}{<}$ &
$\overset{\scriptscriptstyle{\ref{sec:provable}}}{<}$ \\ \hline $\Fin
(\fm)$ & &
$\underset{\rule[-0.5ex]{0pt}{0pt}}{\overset{\rule[2ex]{0pt}{0pt}}{=}}$
& $\overset{\scriptscriptstyle{\ref{sec:custom}}}{>}\
\overset{\scriptscriptstyle{\ref{sec:provable}}}{=}\
\overset{\scriptscriptstyle{\ref{sec:Mostowski}}}{<}\
\overset{\scriptscriptstyle{\ref{sec:Fraenkel}}}{\incomp}$ &
$\overset{\scriptscriptstyle{\ref{sec:custom}}}{>}\
\overset{\scriptscriptstyle{\ref{sec:provable}}}{=}\
\overset{\scriptscriptstyle{\ref{sec:Mostowski}}}{<}\
\overset{\scriptscriptstyle{\ref{sec:Fraenkel}}}{\incomp}$ &
$\overset{\scriptscriptstyle{\ref{sec:provable}}}{<}$ \\ \hline $\Seq
(\fm)$ & & &
$\underset{\rule[-0.5ex]{0pt}{0pt}}{\overset{\rule[2ex]{0pt}{0pt}}{=}}$
& $\overset{\scriptscriptstyle{\ref{sec:provable}}}{=}\
\overset{\scriptscriptstyle{\ref{sec:Mostowski}}}{<}$ &
$\overset{\scriptscriptstyle{\ref{sec:Mostowski}}}{>}\
\overset{\scriptscriptstyle{\ref{sec:provable}}}{\neq}\
\overset{\scriptscriptstyle{\ref{sec:provable}}}{<}\
\overset{\scriptscriptstyle{\ref{sec:Fraenkel}}}{\incomp}$ \\ \hline
$\seq (\fm)$ & & & &
$\underset{\rule[-0.5ex]{0pt}{0pt}}{\overset{\rule[2ex]{0pt}{0pt}}{=}}$
& $\overset{\scriptscriptstyle{\ref{sec:Mostowski}}}{>}\
\overset{\scriptscriptstyle{\ref{sec:provable}}}{\neq}\
\overset{\scriptscriptstyle{\ref{sec:provable}}}{<}\
\overset{\scriptscriptstyle{\ref{sec:Fraenkel}}}{\incomp}$ \\ \hline
$\fzw^\fm$   & & & & &
$\underset{\rule[-0.5ex]{0pt}{0pt}}{\overset{\rule[2ex]{0pt}{0pt}}{=}}$
\\ \hline
\end{tabular}
{}\vspace*{16pt}

One has to read the table from the left to the right and upwards. The
number over a relation refers to the section where the relation was
mentioned.

For any infinite cardinal number $\fm$, if $\Seq (\fm)$, $\seq (\fm)$
and $\fzw^\fm$ are all comparable; the only relations between these
three cardinals which are consistent with ZF are the following:

\begin{enumerate}
\item[(i)] $\Seq (\fm)=\seq (\fm)<\fzw^\fm$\\
(this is true for $\fm=\aleph_0$)
\item[(ii)] $\Seq (\fm)<\seq (\fm)<\fzw^\fm$\\
(see section\;\ref{sec:seq})
\item[(iii)] $\Seq (\fm)<\fzw^\fm<\seq (\fm)$\\
(see section\;\ref{sec:seq})
\item[(iv)] $\fzw^\fm<\Seq(\fm)<\seq (\fm)$\\
(see section\;\ref{sec:Mostowski})
\end{enumerate}

To see this, remember that by Proposition\;\ref{prop:three}\;(2) and
(3), the inequalities $\Seq (\fm)\neq\fzw^\fm$ and $\seq
(\fm)\neq\fzw^\fm$ are both provable in ZF, and further notice that
$\Seq (\fm)=\seq (\fm)$ implies $\aleph_0\le\fm$ which implies
$\fzw^\fm\nleq\Seq (\fm)$ ({\sl cf.}\;\cite[Lemma]{LoriShelah}). So,
in ZF it is provable that there exists no cardinal $\fm$ such that
$\fzw^\fm\le\Seq (\fm)=\seq (\fm)$.

Some other relationships which are provable without the axiom of choice
are the following.

\begin{enumerate}
\item $\fm^2 >\aleph_0\rightarrow \fm > \aleph_0$\\
(see \cite[VIII\,2\;Ex.\,5]{Sierpinski})
\item $\fzw^\fm < \fzw^{\aleph_0}\rightarrow \fm < \aleph_0$ (this means
that $\fm$ is finite)\\
(see \cite[VIII\,2\;Ex.\,3]{Sierpinski})
\item $(\fm\nless\aleph_0\wedge\fm \le \fzw^{\aleph_0})\rightarrow
\fzw^{\aleph_0}\le \fzw^\fm$\\
(see \cite[VIII\,2\;Ex.\,2]{Sierpinski})
\item $\aleph_0 \le \fzw^\fm\rightarrow \fzw^{\aleph_0}\le \fzw^\fm$\\
(see \cite[VIII\,2\;Ex.\,9]{Sierpinski} or Fact\;\ref{fct:aleph0})
\item\label{rel:p36}
$\aleph_0 \le \fzw^\fm\rightarrow \fzw^\fm \nleq \Fin(\fm)^n$
(where $n\in\N$)\\
(see \cite[p.\,36]{LoriShelah})
\item $\aleph_0 \le \fzw^\fm\rightarrow \fzw^\fm \nleq \Fin^n(\fm)$\\
(the proof is similar to the proof of the previous fact \ref{rel:p36})
\item $n\times\Fin(\fm)=\fzw^\fm\rightarrow n=2^k$
(where $n,k\in\N$)\\
(see \cite[p.\,36]{LoriShelah})
\item $\aleph_0 \le \fzw^\fm\rightarrow \fzw^\fm \nleq \Seq (\fm)$\\
(see \cite[Lemma]{LoriShelah})
\item $\aleph_0 \le \fm\rightarrow \fzw^\fm \nleq \seq (\fm)$\\
(the proof is similar to the proof of the Lemma of \cite{LoriShelah})
\item $\fzw^{\fzw^\fm}\neq\fzw^{\aleph_0}$\\
(see \cite[VIII\,2\;Ex.\,7]{Sierpinski})
\item $\big(\fzw^{\Fin(\fm)}\big)^{\aleph_0}=\fzw^{\Fin(\fm)}$\\
(see \cite{Lauchli1})
\item For every $n\in\N$ we have $\aleph_0\nleq\fzw^\fm\rightarrow
n\cdot\fzw^\fm<(n+1)\cdot\fzw^\fm$ \\ (see section\;\ref{sec:misc})
\item $\fzw^{\fzw^\fm}+\fzw^{\fzw^\fm}=\fzw^{\fzw^\fm}$\\
(see \cite{Lauchli1} or section\;\ref{sec:misc})
\end{enumerate}
\newcounter{count}
\setcounter{count}{\value{enumi}}

For each of the following statements we find a permutation model in
which there exists an infinite set $m$ witnessing the corresponding
result, and therefore, by the Jech-Sochor Embedding Theorem, the
following statements are consistent with ZF.

\begin{enumerate}
\addtocounter{enumi}{\value{count}}
\item $n\times \Fin(\fm)=\fzw^{\fm}$
(for any $n\in\N$ of the form $n=2^{k+1}$)\\
(see \cite{LoriShelah})
\item $\aleph_0\le\fzw^{\fzw^\fm}=\fzw^{\Fin(\fm)}$\\
(see \cite[Theorem\;1]{LoriShelah})
\item $\fm^2<[\fm]^2$\\
(see section\;\ref{sec:custom})
\item $\Fin(\fm)<\Seq (\fm)<\fzw^\fm <\seq (\fm)$\\
(see section \ref{sec:seq})
\item $\Seq (\fm)<\seq (\fm)<\Fin(\fm)< \fzw^\fm$\\
(see section \ref{sec:seq})
\end{enumerate}

\nocite{*}

\end{document}